\documentstyle[amssymb,12pt]{article}

\input{tcilatex}

\begin{document}

\begin{center}
\bigskip {\huge On invariant measure structure of }

{\huge a class \ of ergodic discrete dynamical systems.}\bigskip

{\large Anatoliy K. Prykarpatsky}

\bigskip Dept. of Applied Math. at the AGH, Krak\'{o}w 30-059 POLAND,

and Dept. of Nonlinear Math. Analysis

at the IAPMM of the NAS, Lviv 290601, UKRAINA.\\[0pt]
(e-mails: tolik@mailcity.com, prykanat@cybergal.com, prika@mat.agh.edu.pl)\\[%
0pt]
\end{center}

\medskip

\begin{center}
{\footnotesize Abstract. There is studied an invariant \ measure structure
of a class of ergodical discrete dynamical systems by means of the measure
generating function method.}
\end{center}

\section{\protect\large Introduction}

\noindent

\setcounter{equation}{0}Suppose that a topological phase space \ $M$ \ is
endowed with a structure of measurable space \ , that is a $\ \ \sigma -$
algebra \ ${\cal A}(M)$ \ of subsets in \ $M$ , on which there is defined a
finite measure \ $\mu :$\ ${\cal A}(M)$ $\rightarrow {\Bbb {R}}_{+\text{\ }%
}, $ $\mu (M)$ $=1.$ \ As is well known \ [1], \ , a measurable mapping \ $%
\varphi :M$ $\rightarrow M$ \ of the measurable space \ $(M$ , ${\cal A}(M))$
\ \ \ \ is called ergodic \ \ \ discrete dynamical system, if \ $\mu $ -
almost \ everywhere ( $\mu $- a.e.) there exists a not depending on $x\in M$
\ \ limit 
\begin{equation}
\lim_{n\rightarrow \infty }\frac{1}{n}\sum_{k=0}^{n-1}f(\varphi ^{k}x)
\end{equation}
for any bounded \ measurable function \ $f\in {\cal B}(M;{\Bbb {R})}.$

We now assume that the limit (1.1) exists \bigskip\ $\mu $- a.e., that is
one can define a bounded measurable function \ $f_{\varphi }\in {\cal B}(M;%
{\Bbb {R})}$ , where 
\begin{equation}
\lim_{n\rightarrow \infty }\frac{1}{n}\sum_{k=0}^{n-1}f(\varphi
^{k}x):=f_{\varphi }(x)
\end{equation}
for all $\ x\in M$ \ the function (1.2) defines  a finite measure \ $\mu
_{\varphi }:{\cal A}(M)\rightarrow {\Bbb {R}_{+\text{\ }}}$ on $M,$ such
that 
\begin{equation}
\int_{M}f_{\varphi }(x)\text{ }d\mu (x):=\int_{M}f(x)\text{ }d\mu _{\varphi
}(x)\text{.}
\end{equation}
Really, due to the Lebesque-Helley theorem on bounded convergence [2] the
following limit 
\begin{equation}
\int_{M}f_{\varphi }(x)\text{ }d\mu (x)=\lim_{n\rightarrow \infty
}\int_{M}f(x)\text{ }d\mu _{n,\varphi }(x)\text{ ,}
\end{equation}
where by definition,  the Schur averaging 
\begin{equation}
\mu _{n,\varphi }(A):=\frac{1}{n}\sum_{k=0}^{n-1}\mu (\varphi ^{-k}A)\text{ ,%
}
\end{equation}
$n\in {\Bbb {Z}_{+}},$ $A\subset {\cal A}(M),$ and \ $\varphi ^{-k}A:=\{x\in
M:\varphi ^{k}x\in A\}$ for any \ $k\in {\Bbb {Z}_{+}}.$ The limit on the
right hand side of (1.4) exists obviously for any bounded measurable
function \ $f\in {\cal B}(M;{\Bbb {R})}.$ \ \ Thereby \ the equality 
\begin{equation}
\mu _{\varphi }(A):=\lim_{n\rightarrow \infty }\mu _{n,\varphi }(A)\text{ ,}
\end{equation}
for any $\ \ A\subset {\cal A}(M)$ \ defines on the $\ \ \sigma -$ algebra \ 
${\cal A}(M)$ \ an additive not negative bounded mapping $\mu _{\varphi }:%
{\cal A}(M)$ $\rightarrow {\Bbb {R}}_{+}.$ Besides, from the existence \ of
a uniform approximation of arbitrary measurable bounded \ function by means
of finite-valued measurable functions (so called simple ones), one derives
right away the equality (1.3) for any \ $f\in {\cal B}(M;{\Bbb {R})}$. The
needed \ countable additivity of the mapping \ $\mu _{\varphi }:{\cal A}(M)$ 
$\rightarrow {\Bbb {R}_{+\text{\ }}}$ follows from \ the equivalent
expression [6] 
\begin{equation}
\lim_{k\rightarrow \infty }\sup_{n\in {\Bbb {Z}_{+}}}\mu _{n,\varphi
}(A_{k})=0
\end{equation}
for any monotonic sequence of sets \ $A_{j}\supset A_{j+1}$, $j\in {\Bbb {Z}%
_{+}}$ , \ of \ ${\cal A}(M)$ \ with empty transsection.

The measure \ $\mu _{\varphi }:{\cal A}(M)$ $\rightarrow {\Bbb {R}_{+}}$ \
defined by (1.6), has the following invariance property subject to the
dynamical system \ $\varphi :M\rightarrow M:$%
\begin{equation}
\mu _{\varphi }(\varphi ^{-1}A)=\mu _{\varphi }(A)
\end{equation}
for any $\ A\in {\cal A}(M)$ , arising from the following simple identity:

\begin{equation}
\mu _{n,\varphi }(\varphi ^{-1}A)=\frac{n+1}{n}\mu _{n+1,\varphi }(A)-\frac{1%
}{n}\mu (A)\text{ ,}
\end{equation}
after taking the limit  \ \ $n\rightarrow \infty .$ It is easy to see that
(1.8) is completely equivalent to such an equality:

\begin{equation}
\int_{M}f(\varphi x)d\mu _{\varphi }(x)\ =\int_{M}f(x)d\mu _{\varphi }(x)\ \
\ \ \ 
\end{equation}
for any \ $f\in {\cal B}(M;{\Bbb {R})}$. Moreover, if a \ $\sigma -$
measurable set \ $A\in {\cal A}(M)$ is invariant with respect to the mapping
\ \ $\varphi :M\rightarrow M$, that is $\varphi ^{-1}(M)=M$ , then evidently
\ $\mu _{\varphi }(A)=\mu (A)$ . \ 

Therefore the existence of the \ \ $\varphi $-invariant measure \bigskip\ $%
\mu _{\varphi }:{\cal A}(M)$ $\rightarrow {\Bbb {R}_{+}}$ , coinciding with
the measure $\mu :{\cal A}(M)\rightarrow {\Bbb {R}_{+}}$ on the $\sigma -$
algebra \ ${\cal I}(M)$ $\subset {\cal A}(M)$ \ of invariant (with respect
to the dynamical system \ $\varphi :M\rightarrow M$ ) sets, is a necessary
condition of the convergence \ $\mu $- a.e. \ at \ $x\in M$ \ \ of mean
values \ (1.1) as $n\rightarrow \infty $ \ for any $f\in {\cal B}(M;{\Bbb {R}%
).}$ The inverse is also true: due to the Birkhoff theorem [1] , if the
mapping \ $\varphi :M\rightarrow M$ \ conserves \ a finite measure \ $\mu
_{\varphi }:{\cal A}(M)$ $\rightarrow {\Bbb {R}_{+}}$ , mean valueas (1.1)
are convergent \ $\mu _{\varphi }-$ a.e. on \ $M,$ \ and the set of \
convergency is invariant. Thus, if the reduction of the measure \ $\mu :%
{\cal A}(M)\rightarrow {\Bbb {R}_{+}}$ \ upon the invariant \ $\sigma -$
algebra \ ${\cal I}(M)$ $\subset {\cal A}(M)$ \ is absolutely continuous
with respect to that of the measure \ \ $\mu _{\varphi }:{\cal A}(M)$ $%
\rightarrow {\Bbb {R}_{+}}$ , \ then the convergence holds \ $\mu -$ a.e. on
\ $M.$\bigskip

\section{\protect\bigskip {\protect\large A measure generating function
analysis.}}

\bigskip \setcounter{equation}{0}

Assume we are given a discrete dynamical system \ $\varphi :M\rightarrow M$
\ and a sequence of associated measures \ $\mu _{n,\varphi }:{\cal A}(M)$ $%
\rightarrow {\Bbb {R}_{+}}$, $n\in {\Bbb {Z}_{+}},$ defined by (1.5). Then
one can define measure generating functions \ (m.g.f.) $\mu _{n,\varphi
}(\lambda ;A)$, $n\in {\Bbb {Z}_{+}}$, where \ for any \ $A\in {\cal A}(M),$ 
$\lambda \in {\Bbb {C}}$ ,

\begin{equation}
\ \mu _{n,\varphi }(\lambda ;A):=\sum_{k=0}^{n-1}\lambda ^{k}\ \mu (\varphi
^{-k}A).
\end{equation}
Define now the following measure generating \ function

\begin{equation}
\mu _{\varphi }(\lambda ;A):=\lim_{n\rightarrow \infty
}\sum_{k=0}^{n-1}\lambda ^{k}\ \mu (\varphi ^{-k}A),
\end{equation}
where $A\in {\cal A}(M)$ \ and \ $\left| \lambda \right| <1,$ \ that is
needed for finiteness of the expression (2.2). It is easy now to state the
validity of the following lemma.

{\bf Lemma 2.1}{\it . The m.g.f. (2.2) \ satisfies the functional equation\ 
\begin{equation}
\mu _{\varphi }(\lambda ;A)=\lambda \mu _{\varphi }(\lambda ;\varphi
^{-1}A)+\mu (A)
\end{equation}
for any $A\in {\cal A}(M)$ \ and \ $\left| \lambda \right| <1.$ }

{\it $\blacktriangleleft $\bigskip }From (2.3) plainly \ one finds that

\begin{equation}
\mu _{\varphi }(\lambda ;A)-\sum_{k=0}^{n-1}\lambda ^{k}\ \mu (\varphi
^{-k}A)=\lambda ^{n+1}\mu _{\varphi }(\lambda ;\varphi ^{-k-1}A)
\end{equation}
for any $n\in {\Bbb Z}_{+},$ $A\in A(M)$ \ and $\left| \lambda \right| <1.$
Taking the limit in (2.4) as $n\rightarrow \infty ,$\ one arrives at the
expression (2.2) that finishes the proof.$\blacktriangleright $

{\bf Corollary 2.1.} {\it Assume we are given a mapping \ }$\mu _{\varphi
}^{(s)}${\it \ }$:=\mu -s$ $\mu \circ \varphi ^{-1}$ \ {\it on} \ $\ {\cal A}
${\it ${\cal (}$M${\cal )\ }$}, \ {\it where }$\ \left| s\right| <1.$ {\it %
Then the following equality}

\begin{equation}
\mu _{\varphi }^{(s)}(s;A)=\mu (A)
\end{equation}
{\it holds for all }\ \ $A\in $ {\it ${\cal A}$}$(M)$ , $\ \left| s\right|
<1.$

$\blacktriangleleft $ The proof is a straightforward substitution of the
mapping \ $\mu _{\varphi }^{(s)}:{\cal A}${\it ${\cal (}$M${\cal )\
\rightarrow }$}${\Bbb R}$ \ at $\ \left| s\right| <1$ into (2.3).$%
\blacktriangleright $

{\bf Example 2.3.} {\it The induced functional expansion}.

\ \ Let \ $M=[0,1]\subset {\Bbb R}$ \ and \ $\varphi :M\rightarrow M$ \ is
the 'baker'' \ transformation, that is

\begin{equation}
\left\{ 
\begin{array}{c}
\varphi (x):=2x\text{ \ if \ \ \ }x\in \lbrack 0,1/2), \\ 
\varphi (x):=2(1-x)\text{ \ if \ }x\in \lbrack 1/2,1]
\end{array}
\right. .
\end{equation}
Take also another mapping \ $f:M\rightarrow M$ , given as

\begin{equation}
f(x):=2x-x^{2}
\end{equation}
for any $x\in M.$ Then there holds the following decomposition stemming from
(2.5):

\begin{equation}
f(A)=(2-4s)\sum_{n\in {\Bbb Z}_{+}}s^{n}\varphi ^{n}(A)+(4s-1)\sum_{n\in 
{\Bbb Z}_{+}}s^{n}\varphi ^{n}(A)\varphi ^{n}(A)
\end{equation}
for any \ \ $A\in {\cal A}([0,1]).$ In the cases \ $s=1/2$ $\ $and$\ \ s=1/4$
\ one gets simply \ \ for any \ $A\in {\cal A}([0,1])$ \ the decompositions:

\begin{equation}
\sum_{n\in {\Bbb Z}_{+}}(1/2)^{n}\varphi ^{n}(A)\text{ }\varphi
^{n}(A)=f(A)=\sum_{n\in {\Bbb Z}_{+}}(1/4)^{n}\varphi ^{n}(A)\text{ ,}
\end{equation}
which could be useful for some applied set-teoretical considerations. If a
set $\ A=\{x\}\subset \lbrack 0,1]$ \ , from (2.9) one follows such
intersting expansions of the smooth mapping (2.7) into series of piecelike \
linear functions \ (wavelets): 
\begin{equation}
\sum_{n\in {\Bbb Z}_{+}}(1/2)^{n}\varphi ^{n}(x)\text{ }\varphi
^{n}(x)=f(x)=\sum_{n\in {\Bbb Z}_{+}}(1/4)^{n}\varphi ^{n}(x)\text{ .}
\end{equation}
Note here also that a similar expansion given by 
\begin{equation}
\sum_{n\in {\Bbb Z}_{+}}(1/2)^{n}\varphi ^{n}(x)\text{ :=}\xi (x)
\end{equation}
for any $\{x\}\subset \lbrack 0,1]$ , yields the Weierschtrass function \ $%
\xi :[0,1]\rightarrow \lbrack 0,1]$ \ \ nowhere differentiable \ on $%
[0,1]\subset {\Bbb R}$ .

\section{\protect\bigskip {\protect\large The representation of invariant
measure}}

\bigskip \setcounter{equation}{0}3.1. Assume now that there exists the limit
(1.6) \ being due to (1.8) measure preserving on \ ${\cal A}(M)$ . Then the
following theorem is true.

{\bf Theorem 3.1. \ }{\it Let the measure generating function \ }$\mu
_{\varphi }:{\Bbb C}\times $\ ${\cal A}(M)$ $\rightarrow {\Bbb C}$ , {\it %
corresponding \ to a discrete dynamical system \ }$\varphi :M\rightarrow M,$%
{\it \ \ \ exist and satisfy the invariance condition (1.8). Then the limit
expression} 
\begin{equation}
\lim_{
\begin{array}{cc}
\left. \lambda \right\uparrow 1 & (\func{Im}\lambda =0)
\end{array}
}{\it \ }(1-\lambda )\mu _{\varphi }(\lambda ;A)=\mu _{\varphi }(A)
\end{equation}
{\it holds for any} \ $A\in {\cal A}(M)$ .

$\blacktriangleleft \func{Si}$nce all coefficients of the series \ (2.1) \
are bounded, that is of $O(1){\it \bigskip }$ - type, then based on the well
known [3] Hardy theorem, one gets that 
\begin{equation}
\lim_{
\begin{array}{cc}
\left. \lambda \right\uparrow 1 & (\func{Im}\lambda =0)
\end{array}
}{\it \ }(1-\lambda )\mu _{\varphi }(\lambda ;A)=\lim_{n\rightarrow \infty }%
\frac{1}{n}\sum_{k=0}^{n-1}\mu (\varphi ^{-k}A):=\mu _{\varphi }(A)
\end{equation}
for any \ $A\in {\cal A}(M),$ that finishes the proof.$\blacktriangleright $

Therefore we can use the above theorem as the one producing the invariant
measure $\mu _{\varphi }:{\cal A}(M)\rightarrow {\Bbb R}_{+}$ on $M$ \ by
means of the measure generating function \ $\mu _{\varphi }:{\Bbb C}\times $%
\ ${\cal A}(M)$ $\rightarrow {\Bbb C}$ defined by (2.1) for a given discrete
dynamical system \ $\varphi :M\rightarrow M.$ \ It is seen also that the
series (2.1) generates an analytical function as \ \ $\left| \lambda \right|
<1,$ \ having such a property: for any \ $\lambda \in (-1,1)\ \ $and$\ \
A\in {\cal A}(M)$ 
\begin{equation}
\func{Im}\mu _{\varphi }(\lambda ;A)=0.
\end{equation}
Based now on a classical analytical functions theory result \ [4,5] , one
can formulate the following theorem.

{\bf Theorem 3.2.}{\it \ Let a measure generating function} \ $\mu _{\varphi
}:{\Bbb C}\times $\ ${\cal A}(M)$ $\rightarrow {\Bbb C}$ {\it \ satisfy \
the condition (3.3). Then the following representation holds:} 
\begin{equation}
\mu _{\varphi }(\lambda ;A)=\int_{0}^{2\pi }\frac{(1-\lambda ^{2})\text{ }%
d\sigma _{\varphi }(s;A)}{1-2\lambda \cos \text{ }s+\lambda ^{2}}
\end{equation}
{\it for any } $A\in {\cal A}(M),$ {\it where} \ $\sigma _{\varphi }(\circ
;A):[0,2\pi ]\rightarrow {\Bbb R}_{+}$ \ is {\it a function of bounded
variation:} 
\begin{equation}
0\leq \sigma _{\varphi }(s;A)\leq \mu (A)
\end{equation}
{\it for any} $\ s\in \lbrack 0,2\pi ]$ \ {\it and }\ $A\in {\cal A}(M).$

This theorem looks exceptionally interesting for applications since it
reduces the problem of detecting the invariant measure \ $\mu _{\varphi }:$\ 
${\cal A}(M)$ $\rightarrow {\Bbb R}_{+}$ defined \ by 91.6) to calculation
of the following \ complex analytical limit: 
\begin{equation}
\mu _{\varphi }(A)=\lim_{
\begin{array}{cc}
\left. \lambda \right\uparrow 1 & (\func{Im}\lambda =0)
\end{array}
}{\it \ }\int_{0}^{2\pi }\frac{2(1-\lambda )^{2}\text{ }d\sigma _{\varphi
}(s;A)}{1-2\lambda \cos \text{ }s+\lambda ^{2}}\text{ ,}
\end{equation}
where $A\in {\cal A}(M)$ and \ $\sigma _{\varphi }:[0,2\pi ]$ $\times {\cal A%
}(M)\rightarrow {\Bbb R}_{+}$ \ - some Stiltjes measure on $\ \ [0,2\pi ]$,
generated by a given $a$ $priori$ dynamical system \ \ $\varphi
:M\rightarrow M$ and a measure \ $\mu :$\ ${\cal A}(M)$ $\rightarrow {\Bbb R}%
_{+}.$

{\bf Example 3.3}. {\it The Gauss mapping.}

For instance, in the case of the Gauss mapping \ $\varphi :M\rightarrow M,$
where $M=(0,1)$ and for any $x\in (0,1)$ \ \ $\varphi (x):=\{1/x\}$, (here $%
\ "\{\circ \}"$ means taking the fractional part of a number $x\in (0,1)),$
one can show by means of simple but cumbersome calculations, that it is
indeed ergodic \ and possessing the following invariant measure on $M:$%
\begin{equation}
\ \mu _{\varphi }(A)=\frac{1}{\ln 2}\int_{A}\frac{dx}{1+x}\text{ ,}
\end{equation}
yielding obviously, the well known Gauss measure \ $\mu _{\varphi }:$\ $%
{\cal A}(M)$ $\rightarrow {\Bbb R}_{+}$ \ on $M=(0,1).$ As a result, the
following limit for arbitrary \ $f\in L_{1}(0,1)$ is true:

\begin{equation}
\lim_{n\rightarrow \infty }\sum_{k=0}^{n-1}f(\varphi ^{n}x)\stackrel{a.e.}{=}%
\frac{1}{\ln 2}\int_{0}^{1}\frac{f(x)\text{ }dx}{1+x}.
\end{equation}

\section{\protect\large Conclusion}

\setcounter{equation}{0}The analytical expression (3.6) obtained above for
the invariant measure \ $\mu _{\varphi }:$\ $A(M)$ $\rightarrow {\Bbb R}%
_{+}, $ generated by a discrete dynamical system $\varphi :M\rightarrow M,$
looks interesting for concrete calculations. In particular, it is simply
derived from (3.4) that the Stiltjes measure $\sigma _{\varphi }(\circ
;A):[0,2\pi ]\rightarrow \lbrack 0,\mu (A)],$ $A\in {\cal A}(M),$ generates
for any $s\in \lbrack 0,2\pi ]$ a new positive definite measure on \ $A\in 
{\cal A}(M) $ \ as 
\begin{equation}
\sigma _{\varphi }(s)\text{ }(A)=\sigma _{\varphi }(s;A)\text{ , }
\end{equation}
which can be regarded as smearing of the measure \ $\mu :$\ ${\cal A}(M)$ $%
\rightarrow {\Bbb R}_{+}$ along the unit circle ${\Bbb \ S}^{1}$ in the
complex plane ${\Bbb C}$ .

Another important problem still open but related closely with the expression
(3.6), is the following measure evaluation inverse one : how to retrieve the
dynamical system \ $\varphi :M\rightarrow M$ which generated the above
smeared Stiltjes measure \ \ \ \ $\sigma _{\varphi }:[0,2\pi ]$ $\times 
{\cal A}(M)\rightarrow {\Bbb R}_{+}$ \ via the expression \ \ (3.4). Such
and relative problems we plan to discuss in more details elsewhere.

\bigskip

\section{\protect\large Acknowledgments}

\bigskip The author is indebted to Prof. J. Ombach and B. Szafirski from the
Jagiellonian University (Krakow) for discussions. This work was in part
supported through a local \ grant \ of the \ AGH, Krakow.


\bigskip

{\Large References}

\begin{thebibliography}{9}
\bibitem{1}  Wheedon R., Zygmund A. Measure and integral: an introduction to
real analysis. Marcel Decker, Inc., NY,and \ Basel, 1977.

\bibitem{2}  Sinai Ya.G. Ergodic theory. Nauka Publ., Moscow, 1984 (in
Russian)

\bibitem{3}  Hardy G., Convergent series., Cambridge Press, 1947.

\bibitem{4}  Polya G., Sego H. Problems and solutions. Springer, NY, 1982

\bibitem{5}  Privalov I.I. Boundary properties of analytical functions.
Gostekhizdatb Publ., Moscow,1950.

\bibitem{6}  Skorokhod A.V. Elements of the probability theory and casual
processes. Vyshcha Shkola Publ., Kyiv, 1975. \ 
\end{thebibliography}
\end{document}